\documentclass[12pt]{amsart}

\usepackage{bookmark}
\usepackage{graphicx} 
\usepackage{enumerate}
\usepackage{korpi-lib}
\usepackage{tikz-cd}
\usepackage{tikzsymbols}
\usepackage{comment}

\title[Semadeni-Pełczy\'nski derivative and functions on nonmetrizable cubes]{Semadeni-Pełczy\'nski derivative and Banach spaces of \\continuous functions on nonmetrizable cubes}
\author[M. Korpalski]{Maciej Korpalski}
\address{Instytut Matematyczny, Uniwersytet Wroc\l awski, pl.\ Grunwaldzki 2/4, 50-384 Wroc\-\l aw\\ Poland}
\email{Maciej.Korpalski@math.uni.wroc.pl}
\date{}
\thanks{Research supported by Fundación S\'eneca - Agencia de Ciencia y Tecnolog\'ia de la Regi\'on de Murcia (21955/PI/22) and by Młody Badacz grants funded by University of Wrocław}

\begin{document}

\begin{abstract}
We study Banach spaces $C(\prod_{i=1}^n K)$ of real-valued continuous functions on finite products of compact lines. It turns out that the topological character of $K_i$'s translates to an isomorphic invariant of the space of continuous functions on the product.
In particular, for compact lines $K_1, \dotsc, K_n, L_1, \dotsc, L_k$ of uncountable character and $k < n$, we prove that the Banach spaces $C(\prod_{i=1}^n K_i)$ does not embed into $C(\prod_{j=1}^k L_j)$, and there are no continuous linear surjections from $C(\prod_{j=1}^k L_j)$ onto $C(\prod_{i=1}^n K_i)$.

The results are obtained by developing methods used by Semadeni in \cite{ZS60} and Candido in \cite{Ca22}.
\end{abstract}

\maketitle

\section{Introduction} \label{sec:1 Introduction}

A compact line is a linearly ordered set that is compact in the order topology. We study properties of Banach spaces $C(\prod_{i=1}^n K_i)$ of all real-valued continuous functions defined on a finite product of compact lines.
In this paper, we provide a partial answer to the following questions.

\begin{question} \label{question about products} 
Consider compact lines $K_1, \dots, K_n$, $L_1, \dots, L_k$ and $n > k$. 

Are Banach spaces $C(\prod_{i=1}^n K_i)$ and $C(\prod_{j=1}^k L_j)$ isomorphic?

Is there an isomorphic embedding $C(\prod_{i=1}^n K_i) \to C(\prod_{j=1}^k L_j)$?

Is there a continuous linear surjection $C(\prod_{j=1}^k L_j) \to C(\prod_{i=1}^n K_i)$?
\end{question}

Question \ref{question about products} is a generalisation of a famous problem whether the Banach spaces $C([0, 1])$ and $C([0, 1]^2)$ are isomorphic (stated in Banach's book \cite{SB32}). 
We consider mostly nonmetrizable compact lines, as the complete isomorphic characterisation of spaces of continuous functions on metrizable compact spaces has been known since the 1960s, thanks to the results of Miljutin \cite{AM66}, Bessaga and Pełczyński \cite{BP60}.

Let $K_1, \dotsc, K_n, L_1, \dotsc, L_k$ be nonmetrizable compact lines for some $n \neq k$.
Due to the result of Martínez-Cervantes and Plebanek \cite{MCP19}, we know that the products $\prod_{i=1}^n K_i$ and $\prod_{j=1}^k L_j$ are not homeomorphic. 
It follows that the corresponding spaces of continuous functions are not isometrically isomorphic.
In \cite{AM20}, Michalak proved, in particular, that if all the compact lines $K_i, L_j$ are additionally separable, then the Banach spaces $C(\prod_{i=1}^n K_i)$ and $C(\prod_{j=1}^k L_j)$ are not isomorphic.
Thus, it remains to check what happens for nonseparable compact lines. 

In this paper, due to the limitations of our methods, we were able to obtain results only for compact lines of uncountable character.
Nonseparable linearly ordered spaces of countable character can be quite unusual, such as Suslin lines (which are consistently nonexistent in ZFC) or Aronszajn lines (see \cite{Mo09}).

Our approach is inspired by Semadeni's article \cite{ZS60}, which contains a proof that the space of continuous functions on the set of ordinals not greater than the first uncountable ordinal $C([0, \omega_1])$, is not isomorphic to each of its finite powers $\big(C([0, \omega_1])\big)^n$. 
This was one of the first examples of a Banach space that was not isomorphic to its square, addressing another long-standing problem from the book of Banach \cite{SB32}.
Building on top of this result, Candido in \cite{Ca22} defined, for a Banach space $X$, the Semadeni derivative $\cS(X) = X^\mathfrak{s}/X$, where
\[X^\mathfrak{s} = \{x^{**} \in X^{**} \colon x^{**} \text{ is a weak$^*$ sequentially continuous functional}\}.\]
At the end of the Semadeni's paper \cite{ZS60} there is a suggestion by Pełczyński to use a very similar notion, replacing weak$^*$ sequential continuity with continuity on weak$^*$ separable subspaces. 
In this paper, we denote by $\kappa X$ the space of functionals in $X^{**}$ that are weak$^*$ continuous on subspaces of density $\kappa$, and define the $\kappa$-Semadeni-Pełczyński derivative of $X$ (or briefly the SP derivative) as $\cS\cP_\kappa(X) = \kappa X / X$.

The properties of the Semadeni-Pełczyński derivative are a subject of Section \ref{sec: Definitions}. Most of them were proven for the Semadeni derivative by Candido \cite{Ca22}, but several proofs presented here are substantially different.
Later, we introduce the Semadeni-Pełczyński dimension, which measures how many iterations of the derivative are required to obtain the trivial space.

Next, we proceed with some easier results, including the description of the SP derivative of $C(K)$ for a compact line $K$ (Theorem \ref{First applications:5 Characterisation of SP(C(K))}). 
It turns out that the $\kappa$-Semadeni-Pełczyński dimension is nontrivial only for compact lines of character greater than $\kappa$. As a consequence, the character of the compact line $K$ is an isomorphic invariant of the space $C(K)$.
We also present (Theorem \ref{2^theta x [0, lambda+]:4 Galego-like result}) an alternative proof of the results of Galego \cite{Ga09} about the isomorphic classification of spaces $C(2^\theta \times [0, \lambda^+])$ for cardinal numbers $\kappa, \lambda$.

The main result of this paper is a partial answer to Question \ref{question about products}, presented in Section \ref{sec: Calculating SP dim}. 
In Theorem \ref{res:1 Calculating spC(prod K)}, we prove that if $K$ is a finite product of compact lines, then the number of factors of a given uncountable topological character is an isomorphic invariant of the space $C(K)$.
In particular, if compact lines $K_1, K_2, L$ all have an uncountable character, then there are no surjections from $C(L)$ onto $C(K_1 \times K_2)$, nor can $C(K_1 \times K_2)$ be isomorphic to a subspace of $C(L)$.

We finish the paper with Section \ref{sec: Suggestions and remarks}, where we present a collection of thoughts broadly related to the subject of the paper. 
We prove partial results related to Michalak's paper \cite{AM20} on isomorphisms of spaces on continuous functions on products of separable compact lines (Corollary \ref{remarks:4 C(S) not iso to C(Sx prod of compact lines)}), recall some facts about Suslin lines and present a few questions and hypotheses.

The author would like to thank Grzegorz Plebanek for important suggestions concerning the scope of the results in this paper and Antonio Avilés and Gonzalo Martínez-Cervantes for valuable conversations regarding the subject.

\section{Preliminaries} \label{sec:2 Preliminaries}

Every compact space we consider is Hausdorff.
For a compact space $K$, by $C(K)$ we denote the Banach space of all real-valued continuous functions from $K$ considered with the supremum norm.
The dual space $C(K)^*$ can be identified with $M(K)$, the space of signed Radon measures on $K$ with absolute variation as the norm.
In the space $M(K)$, equipped with the weak* topology, the set of Dirac delta measures $\Delta_K = \{\delta_k \colon k\in K\}$ is a topological copy of $K$ (we will sometimes identify $K$ and $\Delta_K$).

If Banach spaces $X$ and $Y$ are isomorphic, we write $X\simeq Y$. If Banach spaces $X, Y$ are isometrically isomorphic, we write $X = Y$, since for our purposes these spaces are identical. For any operator $T : X \to Y$, the formula $T^*y^* (x) = y^*(Tx)$ defines the dual operator $T^* : Y^* \to X^*$. 

For Banach spaces $X, Y$, we write $X \hra Y$ to indicate that $Y$ contains an isomorphic copy of $X$. Similarly, $Y \twoheadrightarrow X$ denotes the existence of a bounded linear surjection from $Y$ onto $X$.

If $\{X_i \colon i \in I\}$ is a family of Banach spaces, we denote by $c_0(I, X_i)$ the $c_0$-product of this family --- the Banach space consisting of all sequences $(x_i)_{i\in I} \in \prod_{i\in I} X_i$ such that for every $\epsilon > 0$, the set $\{i\in I \colon \|x_i\| > \epsilon\}$ is finite, with the supremum norm $\|(x_i)_{i\in I}\| = \sup_{i\in I} \|x_i\|$. 
The dual space of $c_0(I, X_i)^*$ can be seen as $l^1(I, X_i^*)$, and the bidual $c_0(I, X_i)^{**}$ as $l^\infty(I, X_i^{**})$, where the $l^1$- and $l^\infty$-products are defined analogously to the $c_0$-product. A broader description of such spaces can be found in \cite{AA22}*{Section 2, page 5}.

Recall that for a topological space $F$, the character of a point $x \in F$ is the smallest cardinality of a base at the point $x$. The character of a topological space $F$ is the supremum of characters at its points.

For a linearly ordered set $K$, we use the interval notation, so the sets $(a, b), [a, b),$ $ (\leftarrow~\hspace{-4pt}, b], (a, \rightarrow)$ have their usual meanings.
By a compact line, we mean a linearly ordered space with the order topology --- the topology generated by a subbase consisting of all open intervals $(\leftarrow, a), (a, \rightarrow)$.
Recall that all compact lines, as well as their finite products, are sequentially compact. For inequalities between topological coefficients (weight, density, character, etc.) in compact lines, see \cite{EN89}*{Exercise 3.12.4}.

Since ordinal numbers are linearly ordered, we also use the interval notation with respect to them. We write $\omega$ for the set of natural numbers. For any cardinal number $\kappa$, the space $[0, \kappa]$ denotes the ordinal number $\kappa+1$ considered as a compact line. By $\kappa^+$ we mean the smallest cardinal number greater than $\kappa$.

We recall another topological concept introduced by Arhangel’skii in \cite{Ar83}.

\begin{definition}
Let $F$ be a topological space and $\kappa$ a cardinal number. A function $f : F \to \bR$ is called $\kappa$-continuous if its restriction to any subset of $F$ of cardinality at most $\kappa$ is continuous.
\end{definition}

It might be easier to use the following, rather straightforward, characterisation of $\kappa$-continuity.

\begin{lemma}[\cite{Ar83}*{Corollary 1}] \label{Fts:0 Characterisation of k-continuity}
A function $f$ on a topological space $F$ is $\kappa$-continuous if and only if $f|G$ is continuous for every subspace $G\subseteq F$ of density $\le \kappa$.
\end{lemma}

The notion of $\kappa$-continuity is closely related to the functional tightness of a topological space. For any topological space $F$, let $t_0(F)$ be the least cardinal $\kappa$ such that every $\kappa$-continuous function on $F$ is continuous.

It is not very difficult to show that $t_0(F)\le \chi(F)$ for any topological space $F$. In particular, we have $t_0([0, \kappa^+)) = \kappa$ for every infinite cardinal $\kappa$.
In \cite{MK17}, Krupski proved, in particular, that functional tightness is preserved under finite products.
We will use this fact without explicitly referring to it each time.

Now we introduce some definitions that will be used to characterise the space of $\kappa$-continuous functions on products of compact lines. This will allow us to calculate the $\kappa$-Semadeni-Pełczyński derivative of spaces of continuous functions on such products.
We write $A \sqcup B$ for the disjoint union of sets $A$ and $B$.

\begin{definition}
Let $K$ be a compact line and $\kappa$ a cardinal number. A point $k\in K$ is called $\kappa$-inaccessible from the left [right] if $k \notin \ol{A}$ for every set $A\sub (\leftarrow, k)$ [$A\sub (k, \rightarrow)$] of size $\le \kappa$. 

We denote 
\begin{align*}
K_L = \{k \in K &: \chi(k, (\leftarrow, k]) > \kappa\}, \\
K_R = \{k \in K &: \chi(k, [k, \rightarrow)) > \kappa\}, \\
K^\uparrow &= K_L \sqcup K_R,
\end{align*}
so that $K_L$ and $K_R$ are the sets of $\kappa-$inaccessible points of $K$ from the left or right, respectively. The set of all $\kappa$-inaccessible points of $K$ is denoted $K^\uparrow$ (including duplicates, if a point is $\kappa$-inaccessible from both sides). 
Note that the sets $K_L, K_R, K^\uparrow$ all depend on the parameter $\kappa$. 
\end{definition}

The methods developed by Semadeni in \cite{ZS60} have so far been used primarily for ordinal intervals, with the main focus on the space $[0, \omega_1]$. 
The following classical fact motivates the application of results obtained for ordinal intervals to compact lines of uncountable character. 

\begin{fact} \label{Fts:2 Copy of [0, kappa] in K}
Let $\kappa$ be an infinite cardinal number. If $K$ is a compact line and $k \in K^\uparrow$ is a point of character $\geq \kappa$, then we can find a topological copy of $[0, \kappa^+]$ inside $K$ such that $k$ corresponds to the point $\kappa^+$.
\end{fact}

Recall that a topological space $F$ is pseudocompact if every real-valued continuous function on $F$ is bounded, and $F$ is locally compact if it is homeomorphic to an open subset of a compact space.

\begin{theorem} \cite{EN89}*{Theorem 3.10.26} \label{Fts:3 Product of locpsd and psd is psd}
The Cartesian product $F \times G$ of a pseudocompact space $F$ and a pseudocompact locally compact space $G$ is pseudocompact.
\end{theorem}

If $F$ is a Tychonoff space, we denote by $\beta F$ its Stone–Čech compactification. 

\begin{theorem}[Glicksberg] \cite{RW74}*{8.12} \label{Fts:4 Glicksberg products and beta}
If $F$ and $G$ are infinite, then the product space $F \times G$ is pseudocompact if and only $\beta(F \times G) = \beta F \times \beta G$.
\end{theorem}

Let $\kappa$ be any infinite cardinal and $K$ a compact line. 
Due to cofinality reasons, any continuous function $[0, \kappa^+) \to \bR$ must be eventually constant. It follows that $\beta [0, \kappa^+) = [0, \kappa^+]$.

Note that for $k\in K_L$ [$k\in K_R$], the space $(\leftarrow, k)$ [$(k, \rightarrow)$] is a pseudocompact locally compact space. 
By Fact \ref{Fts:2 Copy of [0, kappa] in K}, for any point $k\in K_L$ [$k\in K_R$] we have $\beta (\leftarrow, k) = (\leftarrow, k]$ [$\beta (k, \rightarrow) = [k, \rightarrow)$]. 
Using the two classical theorems above, we claim that a similar property holds for products of compact lines.

\begin{corollary} \label{Fts:5 Beta of the product for k-inaccessible points}
Consider compact lines $K_1, \dotsc, K_n$ and points $k_i \in K_i^\uparrow$ for a given cardinal number $\kappa$. Let 
\begin{equation*}
I_i = 
\begin{cases}
(\leftarrow, k_i)                          &\text{ for } k_i \in K_{i L},\\
(k_i, \rightarrow)                         &\text{ for } k_i \in K_{i R}.
\end{cases}
\end{equation*}
It follows that $\beta \prod_{i=1}^n I_i = \prod_{i=1}^n \beta I_i = \prod_{i=1}^n \overline{I_i}^{K_i}$.
\end{corollary}

\section{New objects and their basic properties} \label{sec: Definitions}

\subsection{Properties of Semadeni-Pełczyński derivative}
First, we define the Semadeni-Pełczyński derivative. Then, we prove its properties, namely that it is preserved by embeddings and continuous surjections, and that it commutes with $c_0$-products.

\begin{definition}
Let $X$ be any Banach space. For an infinite cardinal number $\kappa$ put
\[\kappa X = \{x^{**} \in X^{**} \colon \forall A \sub X^*\ |A| \leq \kappa \ \exists x\in X\ x^{**}|A = x|A\}.\]
The space $\kappa X$ is a closed subspace of $X^{**}$ containing $X$.
We define $\kappa$-Semadeni-Pełczyński derivative by $\cS\cP_\kappa(X) = \kappa X/X$.
\end{definition}

There is some additional context to this definition.
Recall that a topological space $F$ is realcompact if it can be topologically embedded as a closed subspace in the product $\bR^\Gamma$ for some set $\Gamma$.
The space $\omega X$ has already been considered in the literature (see \cite{MT84}*{Section 2.4}), as it is the realcompactification of the Banach space $X$ considered with the weak topology. 
The notation $\omega X$ was introduced by Corson \cite{Co61} (who was actually writing $\aleph_0 X$).

If $Z \sub X$ is a closed subspace of a Banach space $X$, we denote elements of the quotient space $X/Z$ by $[x] = x + Z$ for $x\in X$.
The following standard fact will be crucial for some of the proofs later in this paper.

\begin{lemma} \label{spder:1 Induced operator}
Assume that $X, Y$ are Banach spaces and $T : X \to Y$ is a continuous linear operator. Then the map $\widehat{T} : X/ker(T) \to Y$ given by $\widehat{T}([x]) = T(x)$ is a continuous linear operator of the same norm. If additionally $T$ is surjective, then $\widehat{T}$ is an isomorphism.
\end{lemma}

It turns out that for any compact space $K$ satisfying some natural properties, we can describe elements of the space $\kappa C(K)$ as $\kappa$-continuous functions.
The following theorem is a modified version of \cite{PL91Mazur}*{Proposition, page 29}.

\begin{theorem} \label{Fts:1 Characterisation of kappa C(K)}
Suppose that $K$ is a compact space such that every measure $\mu \in M_1^+(K)$ has a separable set of full measure.
Then $\phi \in \kappa C(K)$ if and only if $\phi$ is represented by a $\kappa$-continuous function on $K$.
\end{theorem}

\begin{proof}
Consider any $\phi \in \kappa C(K)$ and define a function $g_\phi : K \to \bR$ by $g_\phi(k) = \phi(\delta_k)$, for all $k\in K$.
By the definition of $\kappa C(K)$, the function $g_\phi$ is $\kappa$-continuous. 

It remains to show that $\phi(\mu) = \int_K g_\phi \diff \mu$ for all probability measures $\mu\in M(K)$.
By the assumption, for every such $\mu$, there is a closed set $F\subseteq K$ of full measure which is the closure of a countable set $D$.
The desired formula follows from Lemma \ref{Fts:0 Characterisation of k-continuity}.

Conversely, for a $\kappa$-continuous function $f: K \to \bR$, we define $\phi(\mu) = \int_K f \diff \mu$. We can see that $\phi$ is well-defined as $\mu$ has a separable set of full measure.
The functional $\phi$ is an element of $\kappa C(K)$, since the union of supports of $\kappa$ many measures is a set of density $\kappa$ on which, by Lemma \ref{Fts:0 Characterisation of k-continuity}, $f$ is continuous.
\end{proof}

Note that the assumption of Lemma \ref{Fts:1 Characterisation of kappa C(K)} can be relaxed --- it is enough to assume that every probability measure on $K$ has a subspace of density $\kappa$ of full measure. Recall that if $K$ is a compact line, then every measure $\mu\in M_1^+(K)$ has a separable support (see \cite{Me96}*{page 86}). 
This observation easily implies the following.

\begin{lemma}
If $K=K_1\times\dotsc\times K_n$ is a finite product of compact lines $K_i$, then every $\mu \in M_1^+(K)$ admits a separable set of full measure.
\end{lemma}

For the Semadeni derivative, the following properties were proven by Candido (see \cite{Ca22}*{Theorem 3.1, Lemma 3.10 and Theorem 1.2}).
Here, we show that the Semadeni-Pełczyński derivative behaves very similarly to the Semadeni derivative and present slightly different proofs of these facts.
We begin by showing that the SP derivative behaves well with respect to isomorphic embeddings and continuous linear surjections.

\begin{lemma} \label{spder:2 X in Y, SP X in SP Y}
For Banach spaces $X, Y$, if $X \hra Y$, then $\cS\cP_\kappa(X) \hra \cS\cP_\kappa(Y)$ and if $X \twoheadrightarrow Y$, then $\cS\cP_\kappa(X) \twoheadrightarrow \cS\cP_\kappa(Y)$.
\end{lemma}

\begin{proof}
Consider a linear operator $T : X \to Y$. We can lift $T$ to the dual operators $T^{*}: Y^* \to X^*$ and $T^{**} : X^{**} \to Y^{**}$. 
Note that the functionals $\delta_x$ form a closed copy of $X$ in $X^{**}$.
It is a standard to check that $T^{**}x = Tx$ for $x\in X$ and if for some $x^{**} \in X$ we have $T^{**}x^{**} \in Y$, then in fact $x^{**} \in X$.

First, we prove that for $\phi\in \kappa X$ we have $T^{**}\phi\in \kappa Y$. Fix some element $\phi \in \kappa X$ and let $A \sub Y^*$ be any set of cardinality at most $\kappa$. Put $B = T^*[A] \sub X^*$. 
Since the cardinality of $B$ is at most $\kappa$, there exists an element $x\in X$ satisfying $x|B = \phi|B$. It follows that $T^{**}x|A = T^{**}\phi|A$ and $T^{**}\phi \in \kappa Y$. 

Define an operator $\wh{T} : \kappa X \to \cS\cP_\kappa(Y)$ by $\wh{T}(x^{**}) = [T^{**}x^{**}]$, its norm is bounded by $\|T\|$. 
Observe that $X \sub ker \wh{T}$ since for $x\in X$ we have $T^{**}x \in Y$.
By Lemma \ref{spder:1 Induced operator}, we obtain an injective operator $S : \kappa X /ker \wh{T} \to \cS\cP_\kappa(Y)$ of the same norm as $\wh{T}$. 

If $T$ was a surjection, then $S$ is also a surjection, and $\kappa X /ker \wh{T}$ is a quotient of $\cS\cP_\kappa(X)$. Therefore, $\cS\cP_\kappa(X) \twoheadrightarrow \cS\cP_\kappa(Y)$.

Now assume that $T$ is an embedding. If $\wh{T}(x^{**}) = 0$, then $T^{**}x^{**} \in Y$, which implies $x^{**} = x$ for some $x \in X$. It follows that $ker \wh{T} = X$ and $S : \cS\cP_\kappa(X) \to \cS\cP_\kappa(Y)$.

We need to check whether $S^{-1}$ is bounded. Assume towards a contradiction that there exists $\phi \in \kappa X$ such that $\|[\phi]\| = 1$ and $\|\wh{T}\phi\| < \|T^{-1}\|/2$. Then we can pick $y \in Y$ such that $\|T^{**}\phi - y\| < \|T^{-1}\|/2$. 
Note that $dist(T^{**}\phi, T^{**}[X]) \geq \|T^{-1}\|$, as otherwise there is some $x\in X$ such that $\|T^{**}(\phi - x)\| < \|T^{-1}\|$ and $\|\phi - x\| < 1$.
From the triangle inequality it follows that 
\[dist(T^{**}[X], y) \geq dist(T^{**}\phi, T^{**}[X]) - \|T^{**}\phi - y\| > \|T^{-1}\|/2.\]

By the Hahn-Banach theorem, there is an element $y^* \in Y^*$ of norm one such that $y^*(y) > \|T^{-1}\|/2$ and $y^*|T^{**}[X] = 0$. This leads to a contradiction, as
\[\|T^{-1}\|/2 < |y^*(y)| = |\phi(T^*y^*) - y^*(y)| = |(T^{**}\phi - y) (y^*)| \leq \|T^{**}\phi - y\| < \|T^{-1}\|/2.\]
\end{proof}

We now prove that the SP derivative commutes with the $c_0$-product. To do this, we need the following lemma.

\begin{lemma} \label{spder:3 c0 of kappa operation}
If $\{X_i \colon i \in I\}$ is a family of Banach spaces, then for every infinite cardinal number $\kappa$ we have
\[\kappa(c_0(I, X_i)) = c_0(I, \kappa X_i).\]
\end{lemma}

\begin{proof}
Consider an element $(\phi_i)_{i\in I} \in c_0(I, \kappa X_i)$ with only one non-zero coordinate, say $\phi_{i_0} \neq 0$ and $\phi_j = 0$ for $j \in I \sm \{i_0\}$. 
Since every functional from $c_0(I, X_i)^*$ acts on $(\phi_i)_{i\in I}$ in the same way as some functional from $X_{i_0}^*$, we have that $(\phi_i)_{i\in I} \in \kappa(c_0(I, X_i))$ and
\[\|(\phi_i)_{i\in I}\|_{c_0(I, \kappa X_i)} = \|(\phi_i)_{i\in I}\|_{\kappa(c_0(I, X_i))}.\]
It follows that the closed linear span of such vectors must span an isometric copy of $c_0(I, \kappa X_i)$ in $\kappa(c_0(I, X_i))$.

Take any $(\phi_i)_{i\in I} \in \kappa(c_0(I, X_i)) \sub l^\infty(I, X_i^{**})$. 
First, notice that for every $i \in I$ we have $\phi_i \in \kappa X_i$, as each functional from $X_i^*$ can be seen as an element of $c_0(I, X_i)^*$.
Moreover, for every $\eps > 0$, only finitely many elements $\phi_i$ can have norm greater than $\eps$.
Otherwise, we could find a countable sequence of functionals on which no element of $c_0(I, X_i)$ could match the values of $(\phi_i)_{i\in I}$, thus $(\phi_i)_{i\in I} \in c_0(I, \kappa X_i)$.
\end{proof}

\begin{theorem} \label{spder:4 c0 of sp kappa}
If $\{X_i \colon i \in I\}$ is a family of Banach spaces, then for every cardinal number $\kappa$ we have
\[\cS\cP_\kappa(c_0(I, X_i)) = c_0(I, \cS\cP_\kappa(X_i)).\]
\end{theorem}

\begin{proof}
If follows from Lemmas \ref{spder:3 c0 of kappa operation} and \ref{spder:1 Induced operator}.
\end{proof}

\subsection{Semadeni-Pełczyński dimension}
We can define a dimension by iterating the $\kappa$-Semadeni-Pełczyński derivative. Later, we will use it to show the main results of this paper. 
For a Banach space $X$ denote $\cS\cP_\kappa^{(1)}(X) = \cS\cP_\kappa(X)$ and $\cS\cP_\kappa^{(n+1)}(X) = \cS\cP_\kappa(\cS\cP_\kappa^{(n)}(X))$.

\begin{definition} \label{dim:0 definition of the dimension}
We define $\kappa$-Semadeni-Pełczyński dimension of a Banach space $X$ by the following conditions
\begin{itemize}
    \item $sp_\kappa(X) = -1$ if $X = \{0\}$,
    \item $sp_\kappa(X) = n$ if $\cS\cP_\kappa^{(n+1)}(X) = \{0\}$ and $\cS\cP_\kappa^{(n)}(X) \neq \{0\}$,
    \item $sp_\kappa(X) = \infty$ if for all $n \in \omega$ we have $\cS\cP_\kappa^{(n)}(X) \neq \{0\}$.
\end{itemize}
\end{definition}

We have to verify that the $\kappa$-Semadeni-Pełczyński dimension satisfies the properties we require. In particular, we check that it is preserved under isomorphic embeddings and surjections.

\begin{proposition} \label{dim:1 sp and embeddings}
If $X, Y$ are Banach spaces and $X \hra Y$ or $Y \twoheadrightarrow X$, then $sp_\kappa(X) \leq sp_\kappa(Y)$.
\end{proposition}

\begin{proof}
If $X \hra Y$, then by Lemma \ref{spder:2 X in Y, SP X in SP Y}, we have $\cS\cP_\kappa^{(n)}(X) \hra \cS\cP_\kappa^{(n)}(Y)$ for every $n \in \omega$. It follows that if $sp_\kappa(Y) \leq n$, then also $sp_\kappa(X) \leq n$.

If $Y \twoheadrightarrow X$, then again by Lemma \ref{spder:2 X in Y, SP X in SP Y}, we have $\cS\cP_\kappa^{(n)}(Y) \twoheadrightarrow \cS\cP_\kappa^{(n)}(X)$ for every $n \in \omega$. It follows that if $sp_\kappa(X) \geq n$, then also $sp_\kappa(Y) \geq n$.
\end{proof}

The class of Banach spaces for which the $\kappa$-Semadeni-Pełczyński dimension is equal to zero is rather wide. 
In fact, it includes all realcompact Banach spaces and, by Theorem \ref{Fts:1 Characterisation of kappa C(K)}, all spaces $C(K)$ where every measure on the compact space $K$ has a separable support and $t_0(K) \leq \kappa$.
In the sequel, we will also use the fact that this dimension is preserved under $c_0$ products.

\begin{proposition} \label{dim:2 The dimension of c0 products}
Consider a nonempty family of Banach spaces $\{X_i \colon i \in I\}$ and $n\in \omega \cup \{-1\}$.
We have the following.
\begin{enumerate}
    \item If for all $i \in I$ we have $sp_\kappa(X_i) \leq~n$, then we have $sp_\kappa(c_0(I, X_i)) \leq n$.
    \item If there is $i \in I$ such that $sp_\kappa(X_i) \geq~n$, then we have $sp_\kappa(c_0(I, X_i)) \geq n$.
\end{enumerate}
\end{proposition}

\begin{proof}
(1) If $sp_\kappa(X_i) \leq n$ for all $i\in I$, then by Theorem \ref{spder:4 c0 of sp kappa} and Lemma \ref{spder:2 X in Y, SP X in SP Y}, we have
\[\cS\cP_\kappa^{(n)}(c_0(I, X_i)) = c_0(I, \cS\cP_\kappa^{(n)}(X_i)) = \{0\},\]
so $sp_\kappa(c_0(I, X_i)) \leq n$. The proof for the case of (2) is analogous.
\end{proof}

\section{First applications} \label{sec: First applications}

As a warm-up, we calculate the Semadeni-Pełczyński derivative of the space $C(K)$ for a compact line $K$ (Theorem \ref{First applications:5 Characterisation of SP(C(K))}). 
For a cardinal $\theta$, we denote by $2^\theta$ the Cantor cube of appropriate size, not the cardinal exponentiation.
In the latter part of this Section, we use the SP derivative to obtain a partial isomorphic characterisation of spaces $C(2^\theta \times [0, \lambda^+])$ (Theorem \ref{2^theta x [0, lambda+]:4 Galego-like result}). 
This result was inspired by Galego, who in \cite{Ga09} provided (consistently) a complete isomorphic characterisation of these spaces using different methods.

\subsection{Semadeni-Pełczyński derivative and compact lines}
The space of $\kappa$-continuous functions on a compact line $K$ can be described as the space of continuous functions on some other compact line. This new line can be constructed using the $\kappa$-inaccessible points of $K$.

\begin{definition} \label{First applications:1 k-cnt completion}
For a compact line $K$ define its $\kappa$-continuous completion as
\[\bK = K_L \times \{-1\} \cup K \times \{0\} \cup K_R \times \{1\},\]
considered with the lexicographic order and the order topology. It is easy to see that for any $\kappa \geq \omega$, $\bK$ is a compact line of character $\chi(\bK) = \chi(K)$. By $\pi_K$ we denote the projection from $\bK$ to the first coordinate.
\end{definition}

Recall that for a function $f : K \to \bR$, the oscillation at a point $k\in K$ is given by the formula
\[\osc f(k) = \inf \left\{ \sup_{y, z\in V} |f(y) - f(z)| \colon V \ni k \text{ open}\right\}.\]
The oscillation of a function $f$ is the supremum of oscillations at each point, i.e.
\[\osc f = \sup_{k\in K} \osc f(k)\]
The following result links the distance from any bounded function on $\bK$ to the space of continuous functions on $K$ with its oscillation. This connection will be crucial to calculate the norms in the SP derivatives.

\begin{proposition}[\cite{BJ00}*{Proposition 1.18 (ii)}] \label{First applications:2 Distance in l infty}
Let $K$ be a compact space and $f \in l_\infty(K)$. Then the norm of $f$ in the quotient space $l_\infty(K)/C(K)$ is equal to $\frac{1}{2} \osc f$.
\end{proposition}

From Proposition \ref{First applications:2 Distance in l infty} we can easily deduce the following.

\begin{corollary} \label{First applications:3 Extending function from a dense subspace}
Let $A$ be a dense subset of a compact space $K$ and $f : A \to \bR$ be a bounded function. Then $f$ extends to a unique continuous function $\wt{f} : K \to \bR$ if and only if 
\[osc_A f(k) = \inf \left\{ \sup_{y, z\in V \cap A} |f(y) - f(z)| \colon V \ni k \text{ open}\right\} = 0.\]
\end{corollary}

We can now describe the space $\kappa C(K)$ in a more convenient form.

\begin{lemma} \label{First applications:4 Charactersation of kappa C(K)}
For any compact line $K$, there is an isometric isomorphism $T : \kappa C(K) \to C(\bK)$ that satisfies $T[C(K)] = \{f \circ \pi_K: f \in C(K)\}$
\end{lemma}

\begin{proof}
Due to Theorem \ref{Fts:1 Characterisation of kappa C(K)}, we can identify the space $\kappa C(K)$ with the space of real-valued $\kappa$-continuous functions on $K$.

Define an operator $S : C(\bK) \to \kappa C(K)$ by $Sf(k) = f(k, 0)$.

\begin{scclaim}
$Sf$ is a $\kappa$-continuous function.
\end{scclaim}

For every subset $A \sub K$ of cardinality $\kappa$, any point $k \in K_L$ [$k \in K_R$] is not in the closure of the set $A \cap (\leftarrow, k)$ [$A \cap (k, \rightarrow)$]. 
Hence, $Sf$ is $\kappa$-continuous.

We can also see that $S$ is linear and has norm $1$. Moreover, since $K \times \{0\}$ is dense in $\bK$, the operator $S$ is injective.

\begin{scclaim}
The operator $S$ is onto and its inverse $S^{-1}$ has norm $1$.
\end{scclaim}

Fix any $\kappa$-continuous function $g \in \kappa C(K)$.
We show that $g$ extends to a continuous function $\wt{g} \in C(\bK)$.
By Corollary \ref{First applications:3 Extending function from a dense subspace}, it suffices to show that $osc_K g(\Bbbk) = 0$ for every point $\Bbbk \in \bK$.

Consider any point $\Bbbk = (k, b) \in \bK$ and assume that $osc_K g(\Bbbk) > 0$.
If $b = 0$, then either $k \notin K^\uparrow$, $k\in K_L \cap K_R$ or $k \in K_L \triangle K_R$. In all these cases, we have $\chi(\Bbbk) \leq \kappa$, thus $osc_K g(\Bbbk) = 0$, which contradicts our assumption.

Now, without loss of generality, assume that $b = -1$ (case $b = 1$ is symmetric). 
Then we can find two sequences of points $(x_n)_{n\in \omega}, (y_n)_{n\in \omega}$ in $K$ such that for all $n\in \omega$ we have 
\[x_n < y_n < x_{n+1} < k \text{ and } |g(x_n) - g(y_n)| > osc_K g(\Bbbk).\] 
By compactness and $\kappa$-inaccessibility of $\Bbbk$, both sequences $(x_n)_{n\in \omega}, (y_n)_{n\in \omega}$ converge to some point $x\in K$ below $k$, contradicting $\kappa$-continuity of $g$.

It should be clear that $S\wt{g} = g$ and $\|\wt{g}\| = \|g\|$. Note that if $g$ is continuous, then $\wt{g} = g \circ \pi_K$. 
Thus, we may choose $S^{-1}$ as the operator $T$ from the statement.
\end{proof}

It is widely known that if $\bS$ is a double arrow space, then
\[C(\bS)/\{f \circ \pi_{[0, 1]} \colon f \in C([0, 1])\} = c_0(\mathfrak{c}).\]
A similar operation can be performed on any $C(K)$ space where $K$ is a separable compact line. 
This topic was studied in \cite{CSA+21}*{Section 3} and \cite{AM20}, where Michalak called such quotient operators increment-derivative operators and obtained very interesting results by studying their properties.
The following lemma is an analogue of this phenomenon.

\begin{lemma} \label{First applications:5 Characterisation of SP(C(K))}
For any compact line $K$, we have
\[\cS\cP_\kappa(C(K)) = c_0(K^\uparrow).\]
\end{lemma}

\begin{proof}
By Lemma \ref{First applications:4 Charactersation of kappa C(K)}, there is an isometric isomorphism $T : \kappa C(K) \to C(\bK)$, which allows us to identify $\kappa C(K)$ with $C(\bK)$ and $C(K)$ with $T[C(K)] = \{f \circ \pi_K: f \in C(K)\}$.
Put $S : C(\bK) \to c_0(K^\uparrow)$ given for $g \in C(\bK)$ by
\begin{equation*}
Sg(k) =  
\begin{cases}
\left[g(k, 0) - g(k, -1)\right]/2              &\text{ for } k \in K_L, \\
\left[g(k, 1) - g(k, 0)\right]/2               &\text{ for } k \in K_R.
\end{cases}
\end{equation*}

It should be clear that $S$ is a continuous linear operator of the norm $1$. First, we check that $Sg \in c_0(K^\uparrow)$.

Take $g \in C(\bK)$ and assume that for some $\epsilon > 0$ there exists an infinite set $A\sub K^\uparrow$ such that $|Sg(a)| > \epsilon$ for $a \in A$. 
Then $A$ has an accumulation point in $\bK$, at which the oscillation of $g$ cannot be $0$. It follows that $g$ is not $\omega$-continuous --- and therefore not $\kappa$-continuous for $\kappa \geq \omega$.

We can see that for any $g \in C(\bK)$ we have $\|Sg\| = \frac{1}{2}\osc (g|K \times \{0\})$, so $\ker S = C(K)$. 
By Proposition \ref{First applications:2 Distance in l infty}, it follows that $\|[g]\| = \|Sg\|$, where $[g]$ denotes an element of
\[C(\bK)/T[C(K)] = \cS\cP_\kappa(C(K)).\]
Thus, the image of $S$ is closed in $c_0(K^\uparrow)$.

Now, let us show that the image of $S$ contains a dense subset of $c_0(K^\uparrow)$, namely all elements with finite support, denoted $c_{00}(K^\uparrow)$.
Consider the family of intervals
\[\cI = \{((k, -1), \rightarrow) \sub \bK \colon k \in K_L\} \cup \{[(k, 1), \rightarrow) \sub \bK \colon k \in K_R\}.\]
It is standard to check that any element of $c_{00}(K^\uparrow)$ can be represented as $S(\sum_{i=0}^m a_i \chi_{I_i})$ for some $a_i \in \bR$ and $I_i \in \cI$.

Wrapping up everything so far, Lemma \ref{spder:1 Induced operator} gives us a linear isometry
\[\wh{S} : \cS\cP_\kappa(C(K)) \to c_0(K^\uparrow).\]
\end{proof}

In general, the character of the topological space $K$ is not an isomorphic invariant of the space $C(K)$ (the author is aware of an example by Koszmider, which was perhaps not published). However, by Theorem \ref{First applications:5 Characterisation of SP(C(K))}, the situation is different for compact lines.

\begin{corollary}
If $K, L$ are compact lines and $C(K) \simeq C(L)$, then $\chi(K) = \chi(L)$.
\end{corollary}

\subsection{Semadeni-Pełczyński derivative of \texorpdfstring{$C(2^\theta \times [0, \lambda^+])$}{C(2\^theta x [0, lambda+])}}

Let us recall a few classical notions from the theory of large cardinal numbers.
Consider an ultrafilter $\cU$ on some set $A$. Then $\cU$ is $\kappa$-complete for some cardinal $\kappa$ if and only if there is no partition $A = \bigcup_{\alpha < \kappa} X_\alpha$ of $A$ into $\kappa$ many disjoint sets such that $X_\alpha \notin \cU$ for all $\alpha$.

A cardinal number $\kappa$ is measurable if there exists a $\kappa$-complete nonprincipal ultrafilter on $\kappa$ (in other words, a two-valued $\kappa$-additive measure).
We denote the least measurable cardinal by $\fM$ (if it exists). Every measurable cardinal is inaccessible, so their existence is not provable in ZFC. For more details, see \cite{Je03}*{Chapter 10}.

Our interest in measurable cardinals stems from the following result by Talagrand. Note that in this context, by $2^\theta$ we mean the Cantor cube of size $\theta$.

\begin{theorem}[Talagrand, \cite{MT84}] \label{2^theta x [0, lambda+]:1 C(2^omega) is realcompact}
If $\theta < \fM$, then the space $C(2^\theta)$ is realcompact (i.e. $\omega C(2^\theta) = C(2^\theta)$).
\end{theorem}

It is also standard to check that, under the same assumption on $\theta$, we have $\kappa C(2^\theta) = C(2^\theta)$ for every cardinal $\kappa \geq \omega$. 
A simpler proof of Theorem \ref{2^theta x [0, lambda+]:1 C(2^omega) is realcompact} than the original can be found in \cite{Pl91Dyadic}.
It is also a known fact that $t_0(2^\theta) = \omega$ (see, for example, \cite{Pl91Dyadic}*{Theorem 3}).
We now proceed to an analogue of Lemma \ref{Fts:1 Characterisation of kappa C(K)}.

\begin{lemma} \label{2^theta x [0, lambda+]:2 Characterisation of theta C(2^k x [0, l+])}
Consider infinite cardinal numbers $\kappa, \lambda, \theta$. Let $K = 2^\theta \times [0, \lambda^+]$. If $\theta < \fM$, then $\phi \in \kappa C(K)$ if and only if $\phi$ is represented by a $\kappa$-continuous function on $K$.
\end{lemma}

\begin{proof}
Consider any $\phi \in \kappa C(K)$ and define a function $g_\phi : K \to \bR$ by $g_\phi(k) = \phi(\delta_k)$, for all $k\in K$.
By the definition of $\kappa C(K)$, the function $g_\phi$ is $\kappa$-continuous. 
Thus, it remains to check that
\begin{equation}
\phi(\mu) = \int_K g_\phi \diff \mu
\end{equation}
for every measure $\mu \in M(K)$.

For any $\alpha < \lambda^+$, if we write $K_\alpha = 2^\theta\times \{\alpha\}$, then the space $C(K_\alpha)$ embeds into $C(K)$.
It follows from Theorem \ref{2^theta x [0, lambda+]:1 C(2^omega) is realcompact} that $(1)$ holds for measures supported in $K_\alpha$.
Note that, since $[0, \lambda^+]$ is scattered, measures concentrated on $\bigcup_{\alpha \in I} K_\alpha$ for a finite set $I \sub \lambda^+$ form a norm-dense set in $M(K)$, from which we can conclude $(1)$ in the general case.

For the converse, let $g$ be a real-valued $\kappa$-continuous function on $K$. 
The support of every measure $\mu \in M(K)$ is contained in a countable union of $K_\alpha$'s, so for $\kappa$ many measures, $g$ is continuous on their support.
It follows that $\phi$ given by $(1)$ defines an element of $\kappa C(K)$.
\end{proof}

We can now calculate the SP derivative of the space $C(2^\theta \times [0, \lambda^+])$.

\begin{lemma} \label{2^theta x [0, lambda+]:3 Characterisation of SP(C(2^k x [0, l+]))}
Consider infinite cardinal numbers $\theta, \lambda$. If $\theta < \fM$, then
\[\cS\cP_\lambda(C(2^\theta \times [0, \lambda^+])) = C(2^\theta)\]
and
\[\cS\cP_\kappa(C(2^\theta \times [0, \lambda^+])) = 0,\]
for every cardinal $\kappa > \lambda$.
\end{lemma}

\begin{proof}
Consider any $\phi \in \kappa C(2^\theta \times [0, \lambda^+])$ and define a function $g_\phi : 2^\theta \times [0, \lambda^+] \to \bR$ by $g_\phi(x, \alpha) = \phi(\delta_{(x, \alpha)})$, for all $x\in 2^\theta, \alpha < \lambda^+$. 
By Lemma \ref{2^theta x [0, lambda+]:2 Characterisation of theta C(2^k x [0, l+])}, the function $g_\phi$ is $\kappa$-continuous. 
We know that $t_0(2^\theta \times [0, \lambda^+]) = \lambda^+$, so for $\kappa > \lambda$ it follows that $g_\phi$ is continuous.
Thus, we have $\cS\cP_\kappa(C(2^\theta \times [0, \lambda^+])) = 0$. From this point on, assume $\kappa = \lambda$.

We know that $t_0(2^\theta \times [0, \lambda^+)) = \lambda$, so the function $g_\phi^0 = g_\phi|2^\theta \times [0, \lambda^+)$ is continuous. 
By Corollary \ref{Fts:5 Beta of the product for k-inaccessible points}, we have 
\[\beta \left[ 2^\theta \times [0, \lambda^+) \right] = 2^\theta \times [0, \lambda^+],\]
so let $h_\phi$ be the unique continuous extension of $g_\phi^0$ to $2^\theta \times [0, \lambda^+]$. 

Define an operator $S : \lambda C(2^\theta \times [0, \lambda^+]) \to C(2^\theta)$ by the formula
\[S\phi(x) = \left[g_\phi(x, \lambda^+) - h_\phi(x, \lambda^+)\right]/2,\]
for $x \in 2^\theta$. It is immediate that $S$ is continuous and linear. 
Moreover, by Lemma \ref{2^theta x [0, lambda+]:2 Characterisation of theta C(2^k x [0, l+])}, we have $\ker S = C(2^\theta \times [0, \lambda^+])$, as $\phi \in \ker S$ if and only if $g_\phi = h_\phi$, so when $g_\phi$ is continuous.

By Proposition \ref{First applications:2 Distance in l infty}, we have $\|[\phi]\| = \|S\phi\|$ for $[\phi] \in \cS\cP_\lambda(C(2^\theta \times [0, \lambda^+]))$. We can also see that $S$ is a surjection, as for every $f\in C(2^\theta)$ if we put
\begin{equation*}
\wt{f}(x, \alpha) = 
\begin{cases}
2f(x)                       &\alpha = \lambda^+,\\
0                           &\alpha < \lambda^+,
\end{cases}
\end{equation*}
then $S \wt{f} = f$.
By Lemma \ref{spder:1 Induced operator}, there is a linear isometry between $\cS\cP_\lambda(C(2^\theta \times [0, \lambda^+]))$ and $C(2^\theta)$.
\end{proof}

\begin{theorem}  \label{2^theta x [0, lambda+]:4 Galego-like result}
For any infinite cardinal numbers $\theta, \theta', \lambda, \lambda'$, if $\theta$ and $\theta'$ are below the first measurable cardinal, then the spaces $C(2^\theta \times [0, \lambda^+])$ and $C(2^{\theta'} \times [0, \lambda'^+])$ are isomorphic if and only if $\lambda = \lambda'$ and $\theta = \theta'$.
\end{theorem}

\begin{proof}
If $C(2^\theta \times [0, \lambda^+]) \simeq C(2^{\theta'} \times [0, \lambda'^+])$, then due to Lemma \ref{spder:2 X in Y, SP X in SP Y}, we also have
\[\cS\cP_\kappa\Big(C(2^\theta \times [0, \lambda^+])\Big) \simeq \cS\cP_\kappa\Big(C(2^{\theta'} \times [0, \lambda'^+])\Big)\]
for every cardinal $\kappa \geq \omega$. Without loss of generality $\lambda \geq \lambda'$. Now we have two options; in both of them we apply Lemma \ref{2^theta x [0, lambda+]:3 Characterisation of SP(C(2^k x [0, l+]))} with $\theta = \lambda$ to obtain the result. Either
\begin{itemize}
    \item $\lambda > \lambda'$, then $C(2^\theta) \simeq 0$, which is a contradiction, or
    \item $\lambda = \lambda'$, then $C(2^\theta) \simeq C(2^{\theta'})$, which is possible only for $\theta = \theta'$.
\end{itemize}
\end{proof}

\section{Calculating Semadeni-Pełczyński dimension} \label{sec: Calculating SP dim}
It will be convenient to introduce some notation and properly define classes of objects of interest.

\begin{definition} \label{res:0 Definitions of C_k^n, bK}
For $n\in \omega$ define a class $\cC_\kappa^n$ of compact spaces $K$ satisfying the following
\begin{itemize}
    \item $K = \prod_{i=1}^m K_i$ for compact lines $K_1, \dotsc, K_m$ and
    \item $\chi(K_i) > \kappa$ if and only if $i \leq n$.
\end{itemize}
In other words, $\cC_\kappa^n$ is the class of all finite products of compact lines with exactly $n$~factors of character greater than $\kappa$ (maybe up to the permutation of axes). 
For technical reasons, let $\cC_\kappa^{-1}$ be the singleton of the empty set.

If $K \in \cC_\kappa^n$ for some $n$, then denote $\bK = \prod_{i=1}^m \bK_i$, where each $\bK_i$ is the $\kappa$-continuous completion of $K_i$; see Definition \ref{First applications:1 k-cnt completion}. 
Let $\pi_K : \bK \to K$ denote the projection given by $\pi_K((k_i, b_i)_{i=1}^m) = (k_i)_{i=1}^m$.
For $i \leq m$, we will also write $K(i) = \prod_{j = 1, j \neq i}^m K_j$ and $\bK(i) = \prod_{j = 1, j \neq i}^m \bK_j$. 
\end{definition}

Our goal in this section is to prove the following theorem. 

\begin{theorem} \label{res:1 Calculating spC(prod K)}
Consider any compact space $K \in \cC_\kappa^n$ for some $n \in \omega$. Then
\[sp_\kappa(C(K)) = n.\]
\end{theorem}

Note that if $K \in \cC_\kappa^n$ and $L \in \cC_\kappa^m$, then $K \times L \in \cC_\kappa^{n+m}$.
It might be possible to consider, instead of $\cC_\kappa^0$, the class of compact spaces $K$ of functional tightness $\leq \kappa$ with the property that every Radon measure on $K$ has a separable support. However, in this paper, we do not analyse this case. 
Using Proposition \ref{dim:1 sp and embeddings}, we can easily deduce the following.

\begin{corollary} \label{res:2 Applying SP dim to embedings and surjections}
Let $K \in \cC_\kappa^n, L \in \cC_\kappa^m$ for $n > m$. Then
\[C(K) \not \hra C(L) \text{ and } C(L) \not \twoheadrightarrow C(K).\]
\end{corollary} 

In particular, for compact lines of uncountable character $K_1, \dotsc, K_n, L_1, \dotsc, L_m$, if $n > m$, then $C(\prod_{i=1}^n K_i) \not \hra C(\prod_{j=1}^m L_j)$ and $C(\prod_{j=1}^m L_j) \not \twoheadrightarrow C(\prod_{i=1}^n K_i)$.

To prove Theorem \ref{res:1 Calculating spC(prod K)}, we need a few lemmas about the structure of $\kappa C(K)$ and $\cS\cP_\kappa(C(K))$.
We begin with one technical lemma about the extension of functions to $\kappa$-continuous completion. 
Recall that all necessary notation was introduced in Definition \ref{res:0 Definitions of C_k^n, bK}.

\begin{lemma} \label{res:2,5 Extending k-cnt on product K to cnt on product bK}
Consider a product of compact lines $K = K_1 \times \dotsc \times K_n$. 
Then every $\kappa$-continuous function $g : K \to \bR$ has a unique extension to a continuous function on $\wt{g} : \bK \to \bR$. 
\end{lemma}

\begin{proof}
Fix any $\kappa$-continuous function $g : K \to \bR$.
By Corollary \ref{First applications:3 Extending function from a dense subspace}, it suffices to show that $osc_K g(\Bbbk) = 0$ for every point $\Bbbk \in \bK$.
Take any point $\Bbbk = (k_i, b_i)_{i\leq n} \in \bK$ and assume that $osc_K g(\Bbbk) > 0$. 
Denote by $A_\Bbbk$ the set of coordinates $i \leq n$ such that $k_i \in K_i^\uparrow$ and by $B_\Bbbk$ the rest of the coordinates. 
Without loss of generality, assume that $k_i \in K_{iL}$ for every $i \in A_\Bbbk$.

Since every point on a compact line has a descending neighbourhood base, we can construct two sequences of points $(x_\alpha)_{\alpha < \kappa}, (y_\alpha)_{\alpha < \kappa}$ in $K$ such that:
\begin{itemize}
    \item $x_\alpha(i) < y_\beta(i)$ for $i \in A_\Bbbk$ and $\alpha \leq \beta < \kappa$,
    \item $y_\alpha(i) < x_\beta(i)$ for $i \in A_\Bbbk$ and $\alpha < \beta < \kappa$, 
    \item $(x_\alpha(j))_{\alpha < \kappa}, (y_\alpha(j))_{\alpha < \kappa}$ converge both to $x_j = k_j$ for $j \in B_\Bbbk$.
\end{itemize}
Then, for every $i \in A_\Bbbk$, the sequences $(x_\alpha(i))_{\alpha < \kappa}, (y_\alpha(i))_{\alpha < \kappa}$ have a common supremum $x_i \in K_i$ below $k_i$. 
It follows that the function $g$ is not continuous on the set $\{x_\alpha, y_\alpha \colon \alpha < \kappa\} \cup \{(x_i)_{i \leq n}\}$.
\end{proof}

Now we can proceed to the product version of Lemma \ref{First applications:4 Charactersation of kappa C(K)}.

\begin{lemma} \label{res:3 Charactersation of kappa C(prod K)}
For any compact lines $K_1, \dotsc, K_n$, there is an isometric isomorphism $T : \kappa C(K) \to C(\bK)$ that satisfies $T[C(K)] = \{f \circ \pi_K: f \in C(K)\}$.
\end{lemma}

\begin{proof}
Due to Theorem \ref{Fts:1 Characterisation of kappa C(K)}, we can identify the space $\kappa C(K)$ with the space of real-valued $\kappa$-continuous functions on $K$.
Let us denote $\wt{K} = \prod_{i=1}^m (K_i \times \{0\}) \sub \bK$.

Define an operator $S : C(\bK) \to \kappa C(K)$ by $Sf = f|\wt{K}$.

\begin{scclaim}
$Sf$ is a $\kappa$-continuous function.
\end{scclaim}

Let $A \sub K$ be a subset of cardinality $\kappa$. 
For every coordinate $i \leq n$, any point $k \in K_{iL}$ \Big[$k \in K_{iR}$\Big] is not in the closure of the set $\pi_i[A] \cap (\leftarrow, k)$ $\Big[\pi_i[A] \cap (k, \rightarrow)\Big]$,
from which it follows that $Sf$ is $\kappa$-continuous.

We can also see that $S$ is linear, has norm $1$ and since $\wt{K}$ is dense in $\bK$, the operator $S$ is injective.

\begin{scclaim}
The operator $S$ is onto and its inverse $S^{-1}$ has norm $1$.
\end{scclaim}

Fix any $\kappa$-continuous function $g \in \kappa C(K)$.
By Lemma \ref{res:2,5 Extending k-cnt on product K to cnt on product bK}, there is a continuous extension $\wt{g}$ of $g$ on $\bK$.

It should be clear that $S\wt{g} = g$ and $\|\wt{g}\| = \|g\|$. Note that if $g$ is continuous, then $\wt{g} = g \circ \pi_K$. 
Thus, we can choose $S^{-1}$ as the operator $T$ from the statement.
\end{proof}

After establishing a characterisation of $\kappa C(K)$, we proceed to what is perhaps the most difficult result of this paper, from which Theorem \ref{res:1 Calculating spC(prod K)} will follow easily.

\begin{lemma} \label{res:4 Characterisation of SP(C(prod K))}
Let $K \in \cC_\kappa^n$ for some $n\in \omega$. Then the space $\cS\cP_\kappa(C(K))$ is isomorphic to a $c_0$-product of spaces from the class $\cC_\kappa^{n-1}$. 
\end{lemma}

\begin{proof}
By Lemma \ref{res:3 Charactersation of kappa C(prod K)}, there exists an isometry $T : \kappa C(K) \to C(\bK)$ that allows us to identify $\kappa C(K)$ with $C(\bK)$ and $C(K)$ with $T[C(K)] = \{f \circ \pi_K: f \in C(K)\}$. We denote $\wt{K} = \prod_{i=1}^m (K_i \times \{0\}) \sub \bK$.

Let us introduce one more piece of notation. For any $(k, b) \in \bK_i$, let
\[J(k, b) = \prod_{j = 1}^{i-1} \bK_j \times (k, b) \times \prod_{j = i + 1}^m \bK_j.\]
Now define an operator $S : C(\bK) \to \prod_{i=1}^n c_0\Big(K_i^{\uparrow}, C(\bK(i))\Big)$, for $g \in C(\bK)$, by
\begin{equation*}
Sg(i)(k) =
\begin{cases}
\left[g|J_i(k, 0) - g|J_i(k, -1)\right]/2              &\text{ for } k \in K_{i L}, \\
\left[g|J_i(k, 1) - g|J_i(k, 0)\right]/2               &\text{ for } k \in K_{i R}.
\end{cases}
\end{equation*}
It is standard to check that $S$ is a continuous linear operator. We will prove that $S$ is well-defined and $\ker S = C(K)$.

Fix $i \leq n$ and any function $g \in C(\bK)$. Suppose that for some $\epsilon > 0$, there is a countable infinite set $A \sub K_i^{\uparrow}$ such that $\|Sg(i)(a)\| > \epsilon$ for $a \in A$. 
Then $A$ has an accumulation point in $K_i$ at which the oscillation of $g$ cannot be $0$, so $g$ is not $\omega$-continuous (and also $\kappa$-continuous for $\kappa \geq \omega$).

Take any $g \in C(\bK)$ and let 
\[H = \{(k_i)_{i \leq n} \in K \colon \exists i \leq n \ k_i \in K_i^\uparrow\}\]
be the set of points in $K$ with at least one $\kappa$-inaccessible coordinate. It should be clear that 
\[\osc g|\wt{K} = \sup_{k \in H} \max_{\Bbbk_1, \Bbbk_2 \in \bK} \{f(\Bbbk_1) - f(\Bbbk_2) \colon \pi_K (\Bbbk_1) = \pi_K (\Bbbk_2) = k\},\]
as outside the copy of $H$ in $\wt{K}$, the function $g|\wt{K}$ is continuous. 
Consider any $k\in H$. For some $d \leq 2n$ we can find a sequence of points $\Bbbk_{j\leq d}$ in $\bK$ such that, for each $j \leq d$, the value $|f(\Bbbk_j) - f(\Bbbk_{j+1})|/2$ is, up to a sign, obtained by $Sg$ and the oscillation of $g|\wt{K}$ at $k$ is equal to $f(\Bbbk_1) - f(\Bbbk_d)$. 
It follows that $\|Sg\| \geq \frac{1}{4n}\osc g|\wt{K}$, so $\ker S$ contains only those functions $g \in C(\bK)$ for which the oscillation of $g|\wt{K}$ is $0$. It is also clear that for any $g \in C(K)$, we have $Sg = 0$, thus $\ker S = C(K)$.

By Proposition \ref{First applications:2 Distance in l infty}, if for $g \in C(\bK)$ we consider $[g]$ as an element of
\[C(\bK)/T[C(K)] = \cS\cP_\kappa(C(K)),\]
then $\|[g]\| = \frac{1}{2} \osc g|\wt{K}$. We have already shown that 
\[\|Sg\| \geq \frac{1}{4n} \osc g|\wt{K} = \frac{1}{2n}\|[g]\|,\]
which implies that the image of $S$ is closed in $c_0(K^\uparrow, C(\bK(i)))$.

Now we will show that the image of $S$ is dense in its codomain. By linearity, it is enough to show that for each $i \leq n$, the image of $S$ is dense in $c_0(K_i^\uparrow, C(\bK(i)))$.
Our goal is to show that all elements with only finitely many non-zero coordinates, denoted $c_{00}(K_i^\uparrow, C(\bK(i)))$, lie in the image of $S$.
For $k \in K_L$ consider a clopen rectangle 
\[P^i_k = \{(k_i, b_i)_{i\leq n} \in \bK \colon (k_i, b_i) > (k, -1)\}\]
and for $k\in K_R$ a clopen rectangle
\[R^i_k = \{(k_i, b_i)_{i\leq n} \in \bK \colon (k_i, b_i) > (k, 0)\},\]
then by $\cI_i$ denote the family of all these rectangles, i.e.
\[\cI_i = \{P^i_k \colon k \in K_L\} \cup \{R^i_k \colon k \in K_R\}.\]
It is standard to verify that any element of $c_{00}(K_i^\uparrow, C(\bK(i)))$ can be represented as $S(\sum_{j=0}^m f_j \chi_{I_j})$, for some $I_j \in \cI_i$ and $f_j \in C(\bK(i))$. 

Wrapping up everything so far, Lemma \ref{spder:1 Induced operator} gives us an isomorphism
\[\wh{S} : \cS\cP_\kappa(C(K)) \to \prod_{i=1}^n c_0(K_i^\uparrow, C(\bK(i))).\]
\end{proof}

We now have all the tools necessary to calculate the Semadeni-Pełczyński dimension for spaces $C(K)$, where $K$ is a product of compact lines. 

\begin{proof} (of Theorem \ref{res:1 Calculating spC(prod K)})
We will prove the statement by induction on $n$. 

By Theorem \ref{Fts:1 Characterisation of kappa C(K)}, for $K \in \cC_\kappa^0$, we have $\kappa C(K) = C(K)$, so $sp_\kappa(C(K)) = 0$.

Now assume that for some $n \geq 1$ and every $K \in \cC_\kappa^{n-1}$ we have $sp_\kappa(C(K)) = n - 1$. 
Consider any space $K \in \cC_\kappa^n$. Then there are some compact lines $K_1, \dotsc, K_m$ such that $K = \prod_{i=1}^m K_i$, and $\chi(K_i) > \kappa$ if and only if $i \leq n$. 

By Lemma \ref{res:4 Characterisation of SP(C(prod K))}, we have $\cS\cP_\kappa(C(K)) \simeq c_0(I, C(L_i))$ for some set $I$ and spaces $L_i \in \cC_\kappa^{n-1}$.
From Proposition \ref{dim:2 The dimension of c0 products} it follows that
\[sp_\kappa \left[ c_0(I, C(L_i))\right] = \sup\{sp_\kappa\left(C(L_i)\right) \colon i \in I\}.\]
By the induction hypothesis $sp_\kappa\left(C(L_i)\right) = n - 1$ and thus
\[sp_\kappa\left(C(K)\right) = 1 + sp_\kappa(\cS\cP_\kappa(C(K))) = 
1 + \sup\{sp_\kappa\left(C(L_i)\right) \colon i \in I\} = 1 + n - 1 = n.\]
\end{proof}

\section{Remarks and problems} \label{sec: Suggestions and remarks}

In this Section, we present a collection of scattered thoughts, remarks and questions broadly related to the subject of the paper.
All results follow from already known facts, so the proofs should be straightforward and may get quite sketchy.
Moreover, some of the remarks made here might be known in the broader community.

\subsection{Remarks on the double arrow space}
Denote the double arrow space by $\bS$. Recall that $\bS$ can be represented as the lexicographic product of the unit interval and the two-point space, i.e., $[0, 1] \times_{lex} \{0, 1\}$ (typically excluding the two isolated endpoints).
In what follows, we present a few Remarks on isomorphisms between spaces of continuous functions on $\bS$ and its products.

Regarding the next proposition, if one is familiar with the concept of free dimension defined in \cite{MCP19}, it seems worth noting that free-dim$(\bS) = 1$, while free-dim$(\bS \times [0, \omega]) = 2$.
From Proposition \ref{remarks:1 C(S) iso to C(Sx[0,omega])} it follows that the free dimension of a compact space $K$ is not an isomorphic invariant of the Banach space $C(K)$.

\begin{proposition} \label{remarks:1 C(S) iso to C(Sx[0,omega])}
$C(\bS) \simeq C(\bS \times [0, \omega]).$
\end{proposition}

\begin{proof}
The space $\bS$ contains nontrivial convergent sequences, which implies that the space $C(\bS)$ is isomorphic to its hyperspaces. Denote $I_n = [(1/(n+1), 1), (1/n, 0)] \sub \bS$ and let $\phi_n$ be the homeomorphism between $\bS$ and $I_n$. Then the formula
\[Tf = ((f|I_n) \circ \phi_n)_{n\in \omega}\]
defines an isometry between a hyperplane $C'(\bS) = \{f \in C(\bS) \colon f((0, 1)) = 0\}$ and the space $c_0(\omega, C(\bS))$. It follows that
\[C(\bS) \simeq C'(\bS) \simeq c_0(\omega, C(\bS)) \simeq C(\bS \times [0, \omega]).\]
\end{proof}

By Bessaga and Pełczyński's characterisation of isomorphisms of spaces of continuous functions on countable compacta \cite{BP60}, it follows that $C([0, \omega]) \not \simeq C([0, \alpha])$ for any ordinal $\omega^\omega \leq \alpha < \omega_1$.
This raises the following question.

\begin{question}
Is the space $C(\bS)$ isomorphic to $C(\bS \times [0, \alpha])$ for $\omega^\omega \leq \alpha < \omega_1$?
\end{question}

However, it is interesting to note that Proposition \ref{remarks:1 C(S) iso to C(Sx[0,omega])} does not hold if $[0, \omega]$ is replaced by any nonscattered compact space.
Note that the proof of the following result is similar to that of Lemma \ref{res:4 Characterisation of SP(C(prod K))}.

\begin{proposition} \label{remarks:2 C(S) not iso to C(SxK)}
If $K$ is a compact space such that $C(K)$ does not embed into $c_0(\mathfrak{c})$, then $C(\bS) \not \simeq C(\bS\times K).$
\end{proposition}

\begin{proof}
Assume that $T : C(\bS) \to C(\bS\times K)$ is an isomorphism. 
There is a natural copy of $C([0, 1])$ embedded in $C(\bS)$ (of functions constant on the second coordinate).
Let $Y = T[C([0, 1])]$; this is a separable subspace of $C(\bS\times K)$. 
It is well known that the quotient $C(\bS) / C([0, 1])$ is isomorphic to $c_0(\mathfrak{c})$; see \cite{CSA+21}*{Section 3}.

Let us show that there exists an isomorphic copy of $C(K)$ inside $C(\bS\times K) / Y$. 
Denote by $D$ a countable dense subset of $Y$. For each $f\in C(\bS\times K)$, define the section function $f_{s}(k) = f((s, k))$ for $s\in \bS, k\in K$ and put
\[A_x = \{f\in D \colon f_{(x, 0)} \neq f_{(x, 1)}\}\]
for $x\in [0, 1]$.

Since $D$ is countable, we can fix some $x\in [0, 1]$ such that $A_x = \emptyset$. Then for all $f\in Y$, we have $f_{(x, 0)} = f_{(x, 1)}$. Define
\[\wt{f}(s, k) = \chi_{[(x, 1), \rightarrow)}(s) \cdot f(k)\]
for $f \in C(K), s\in \bS$ and $k\in K$.
Then the subspace $C_x = \{\wt{f} \colon f \in C(K)\}$ defines an isomorphic copy of $C(K)$ inside $C(\bS\times K) / Y$, as for every $f\in C_x$ we have $dist(f, Y) \geq \|f\|/2$.
This contradicts the assumption that $C(K) \not \hra c_0(\mathfrak{c})$.
\end{proof}

Let $K$ be any compact space.
Recall that $C(K)$ is an Asplund space if and only if $K$ is scattered, and that every closed subspace of an Asplund space is Asplund.
It follows from Proposition \ref{remarks:2 C(S) not iso to C(SxK)} that for any nonscattered compact space $K$ we have $C(\bS) \not\simeq C(\bS\times K)$, since $c_0(\mathfrak{c})$ is an Asplund space.

In particular, $K$ can be any product of uncountable compact lines or the unit interval. Thus, we have obtained an easy way to see some of the results of Michalak from \cite{AM20}.

\begin{corollary} \label{remarks:3 C(S) not iso to C(Sx[0, 1])}
$C(\bS) \not \simeq C(\bS\times [0, 1])$.
\end{corollary}

\begin{corollary} \label{remarks:4 C(S) not iso to C(Sx prod of compact lines)}
If $K$ is any product of nonmetrizable separable compact lines, then $C(\bS) \not \simeq C(\bS\times K)$.
\end{corollary}

A closer inspection of the proof of Proposition \ref{remarks:2 C(S) not iso to C(SxK)} shows that, instead of $\bS$, we can put any separable compact line of weight greater than $\omega$ to obtain a similar result.

\subsection{Remarks on the Suslin lines}
Using classical results, we can prove one simple observation about the isomorphic structure of spaces of continuous functions on the products of Suslin lines. The first thing to note is the following classical result by Kurepa.

\begin{theorem}[Kurepa, \cite{Ku50}]
The square of a Suslin line is not ccc.
\end{theorem}

More details on chain conditions can be found in the survey by Todorčević \cite{To00}. It is worth noting that the ccc property of a compact space translates to an isomorphic property of its space of continuous functions.

\begin{theorem}[Rosenthal \cite{Ro69}]
A compact space $K$ is ccc if and only if every weakly compact subset of $C(K)$ is separable.
\end{theorem}

Combining the theorems of Kurepa and Rosenthal, we easily get the following.

\begin{corollary}
Let $\cS$ be a Suslin line. Then $C(\cS) \not\simeq C(\cS^n)$.
\end{corollary}

It is worth mentioning that the above reasoning does not need to hold for the product of two distinct Suslin lines, as such a product can be ccc, according to the results of Jensen (as claimed by, e.g., Rudin in \cite{Ru78}).

\subsection{Remark on the Semadeni-Pełczyński dimension}
It is not clear whether Definition \ref{dim:0 definition of the dimension} can be extended to infinite ordinals. It should probably be done using either the direct or inverse product of Banach spaces, but it is not clear how to define the corresponding mappings. Moreover, it seems reasonable that the space
\[c_0(\omega, C([0, \omega_1]^n))\]
or a very similar space considered by Candido in \cite{Ca22}*{Theorem 1.4} should have the $\omega$-Semadeni-Pełczyński dimension equal to $\omega$.

\bibliographystyle{amsalpha}
\bibliography{refs}

\end{document}